\documentclass[amssymb]{amsart}
\usepackage{amssymb}
\usepackage{latexsym}
\usepackage{xypic}
\makeatletter
\renewcommand{\@makefnmark}{}
\makeatother
\title[Kontsevich's conjecture]{On Kontsevich's Hochschild cohomology 
conjecture}
\author[P. Hu, I. Kriz, A.A. Voronov]{P. Hu, I. Kriz and A.A. Voronov}
\input xypic
\newtheorem{theorem}{Theorem}
\newtheorem{proposition}[theorem]{Proposition}

\newtheorem{lemma}[theorem]{Lemma}

\def\Proof{\medskip\noindent{\bf Proof: }}

\def\umu{\underline{\mu}}

\def\Z{\mathbb{Z}}

\def\R{\mathbb{R}}

\def\Pi{\mathbb{P}^{\infty}}

\def\qed{\hfill$\square$\medskip}

\def\Zpk{\mathbb{Z}/p^{k}}
\def\Zpk1{\mathbb{Z}/p^{k-1}}

\def\mcpk{\mathcal{C}^{\prime}_{k}}
\newcommand{\rref}[1]{(\ref{#1})}

\newcommand{\cform}[3]{\begin{array}{c}
{\scriptstyle #3}\\
#1\\
{\scriptstyle #2}\end{array}}
\def\r{\rightarrow}
\def\mc{\mathcal{C}}
\def\omc{\overline{\mathcal{C}}}
\def\mb{\mathcal{B}}

\def\sl2{\widetilde{SL_{2}(\Z)}}

\begin{document}
\maketitle


\vspace{5mm}

\section{Introduction}

\footnote{The authors were supported by the NSF.}
A conjecture of Deligne stated that the Hochschild cohomology complex
of an 
associative algebra has a natural structure of a $2$-algebra, i.e.
an algebra over the chain complex version of
the $2$-cube operad. 
This indicated a remarkable connection between the deformation
theory of associative algebras, and the geometry of configuration spaces
of points in the plane. There are several known proofs of Deligne's conjecture,
see \cite{tam2}, \cite{tam3}, \cite{vor},
\cite{ms}, \cite{kon1}. The purpose of the present paper is to prove a 
generalization conjectured by
Kontsevich \cite{kon}, calling
for an analogue of Deligne's conjecture for algebras over the
little $k$-cube operad.  

\vspace{3mm}
The first problem is to define a suitable generalization of the 
Hochschild cohomology complex. Kontsevich \cite{kon}
proposes to do this by modifying the Quillen cohomology
complex, but that approach forces some restrictions 
(in fact, it only seems to work for the little $k$-cube operad, 
and, as stated, only in chain complexes 
over fields of characteristic $0$). 
A key feature of our approach is that we give  
a completely natural definition of the Hochschild cohomology complex,
not restricted to those situations.
In fact,
all of our constructions in principle work for any operad, and in
a closed symmetric monoidal 
category (\cite{ml}).
Nevertheless, 
to avoid technical problems, we shall still stick to specific
cases. Namely, in the statement of our theorems,
we shall assume that 
{\em $\mb$ is the category of sets, or $K$-modules where $K$ is a field}
(not necessarily of characteristic $0$).
There is at least one substantially different case of interest, namely
the case of spectra ($S$-modules \cite{ekmm}). However, 
homotopical algebra in that case is more difficult, and
will not be discussed here.

\vspace{3mm}
Now consider operads $\mc$ in the category of simplicial sets. An 
example of special interest to us is the operad $\mc_{k}$ which is the {\em
set of singular simplices of the operad
of little $k$-cubes} \cite{271}. 
In the beginning of the next Section, we will introduce the notion
of $\mc$-algebras $R$ in the category $s\mb$, and $(\mc,R)$-modules.
Furthermore, for $(\mc,R)$-modules $M,N$, we will construct a
`derived mapping
object'  
$$RHom_{(\mc,R)}(M,N)$$
in the homotopy category of $s\mb$. Then our main result can be stated
as follows:

\begin{theorem}
\label{t1}
Let $R$ be a cofibrant $\mc_{k}$-algebra (in $s\mb$ as above). 
Then there is a functorial model of
$$RHom_{(\mc_{k},R)}(R,R)$$
which is a 
$\mc_{k+1}$-algebra.
\end{theorem}

\vspace{3mm}
Tamarkin \cite{tam1} 
previously proved a purely algebraic version of Theorem \ref{t1}
for the Quillen cohomology complex
\cite{q} (also known as the deformation complex), 
working in the category of chain
complexes rather than simplicial modules.
His proof works over fields of characteristic 
$0$, and does not include the statement  of Deligne's original
conjecture (the case of $k=1$).
It uses the fact that the little cube chain operad in that case is
formal and Koszul. 

\vspace{3mm}
Theorem \ref{t1} is based on the geometry of the
little $k$-cube operad. Its proof will be
completely `derived' and based on a general principle
(cf. \cite{dunn}) that ``commuting $\mathcal{C}_k$- and $\mathcal{C}_1$-structures
give a $\mathcal{C}_{k+1}$-structure''. In the case of Hochschild cohomology
of $R$, the $\mathcal{C}_{k}$-structure is induced by the 
$\mathcal{C}_k$-structure of $R$, the $\mathcal{C}_1$-structure
from the Yoneda product. Our proof applies to modules over fields of
any characteristic (as well as sets), and does include
the case $k=1$. In the statement of the Theorem, the assumption that
$R$ be cofibrant should not be regarded as a restriction, since for
non-cofibrant algebras the correct notion of Hochschild cohomology is
obtained by first taking a cofibrant replacement: this leads to 
the right notion of Hochschild cohomology in the case $k=1$.

\vspace{3mm}
The present paper is organized as follows: In the next Section, we shall
reduce Theorem \ref{t1} to a much more general context, which may be of
independent interest as a generalization of the Kontsevich conjecture. 
We will introduce the notion of operads fibered over
a given operad $\mc$ in simplicial sets. Using
Boardman-Vogt's
$\Box$-product of operads \cite{BV}, 
we will state in Theorem \ref{t2} that
for a certain class of operads $Q$ fibered over $\mc$, which we call 
{\em special},
there is a notion of Hochschild cohomology object constructed from $Q$,
and that moreover
this object has the structure of a $\mc\Box\mc_{1}$-algebra.
We also prove a simplicial version of the theorem of G. Dunn \cite{dunn}
comparing $\mc_k\Box\mc_{\ell}$ to $\mc_{k+\ell}$.

\vspace{3mm}
There are three facts which together will reduce the proof of Theorem \ref{t1}
to Theorem \ref{t2}: First, a $\mc_{k}$-algebra gives rise to an example
of a special operad fibered over $\mc_{k}$. Second, the two ensuing notions
of Hochschild cohomology coincide. Third, we have $\mc_{k}\Box\mc_{1}\simeq
\mc_{k+1}$. All these statements are proven in Section \ref{s2} with the
exception of the `special' property: that is technical, and left to
Section \ref{sspec}.

\vspace{3mm}
In Sections \ref{s3}, \ref{s4}, we will describe the technical machinery used
to prove Theorem \ref{t2}. This technique can be separated into two steps: 
First,
in Section \ref{s3}, we shall consider ``lax algebras'' over an operad $\mc$
and show how they may be turned into strict algebras over a different,
but weakly equivalent, operad $\overline{\mc}$. In Section \ref{s4},
we shall introduce an additional, ``vertical'' category structure on
a lax algebra, which will allow us to get a 
$\overline{\mc}\Box\mc_{1}$-algebra.
Both constructions must be discussed also in categories ``enhanced'' over
a given category $\mb$ (we deliberately choose a different
name than ``enriched'' which is a different
notion, although also related
to our considerations). In Section \ref{s4a}, we will apply the techniques
of Sections \ref{s3}, \ref{s4} to our main example, which will give the
proof of Theorem \ref{t2}.

\vspace{3mm}
\noindent
{\bf Acknowledgement:} We are very indebted to T. Fiore for comments.

\vspace{5mm}
\section{Special operads fibered over $\mathcal{C}$}
\label{s2}

In this Section, we shall reduce Theorem \ref{t1} to another statement,
which can be phrased in a more general context. However, we begin by
filling in the missing definitions in the statement of Theorem \ref{t1}.

Recall \cite{ml} that a {\em closed symmetric monoidal category}
is a category $\mathcal{B}$ with a symmetric monoidal structure
$\boxtimes$ and an `internal Hom'
functor $Hom:\mb^{Op}\times\mb\r\mb$ with a natural
bijection
\begin{equation}
\label{ecl}
Mor(A\boxtimes B,C)\cong Mor(A, Hom(B,C)),
\end{equation}
satisfying certain axioms.
We shall work in the category $s\mb$
of simplicial objects in $\mb$.

For a set $S$ and an object (or morphism) $X$ of $\mb$, let
$$S\otimes X = \cform{\coprod}{S}{} X.$$
We say that an object $R$ of $\mb$ is a $\mc$-algebra for an operad
$\mc$ in sets if there are structure maps
\begin{equation}
\label{e1}
\mc(n)\otimes(R\boxtimes...\boxtimes R)\r R
\end{equation}
satisfying the usual axioms (see e.g. \cite{km}). Now an important point
is that here (and analogously in other places below),
the same definition may be used for an operad in simplicial
sets and an object of $s\mb$ (this simply means that
maps \rref{e1} exist on each simplicial level, and are natural
with respect to simplicial structure). For a $\mc$-algebra $R$,
a $(\mc,R)$-module is an object $M$ of $\mb$ together with structure maps
\begin{equation}
\label{e2}
\mc(n)\otimes(M\boxtimes R\boxtimes...\boxtimes R)\r M
\end{equation}
satisfying the usual axioms (c.f. \cite{km}). 
If $C_{R}$ is the monad in $s\mb$
defining free $(\mc,R)$-modules, then define for $(\mc,R)$-modules $M,N$,
$$Hom_{(\mc,R)}(M,N)$$
as the equalizer of the two obvious maps
$$Hom(M,N)\begin{array}{c}\r\\ \r\end{array}Hom(C_{R}M,N).$$
(One is induced by $C_{R}M\r M$, the other by the natural map
$Hom(M,N)\r Hom(C_{R}M,C_{R}N)$, composed with the map induced by 
the map $C_{R}N\r
N$.) 

\vspace{3mm}
We will be considering Quillen (closed) model structures on certain categories, which
will be needed to do homotopy theory in those categories. A (closed) model
structure on a complete cocomplete
category consists of three classes of morphisms called
fibrations, cofibrations and equivalences. (The word closed may be omitted
as it carries no meaning.) All the information
on closed model structures needed in this paper, and all the methods
of construting closed model structures we shall need can be found in 
\cite{qh}, \cite{ds}. One should remark that on many categories there
are many closed model structures which have the same equivalences
and hence lead to the same homotopy theory. However, when doing
constructions, one typically needs to fix a closed model structure in
order to control the homotopy behavior of the objects produced by
the construction.

\vspace{3mm}
The most basic case needed is the category of simplicial
sets (see \cite{ds}, 11.1). In this category, there is a closed model structure
where cofibrations are injective maps (more precisely sequences
of injective maps), and equivalences are maps of simplicial sets
which produce homotopy equivalences after simplicial realization.
(It then follows that fibrations are so called Kan fibrations,
but that is less important to us.)
Next, recall (\cite{ds}, 11.2) that a Quillen model structure 
in $s\mb$ is obtained as follows:
Let $U$ be the forgetful functor from $s\mb$ to simplicial sets. Then
a morphism $f$ in $s\mb$ is a fibration or equivalence if and only if
$Uf$ is a fibration (equivalence). We say that fibrations and equivalences
are {\em created} by the forgetful functor $U$.

\vspace{3mm}
\noindent
{\bf Remark:}
Actually, it turns out that in the two
cases we consider in this paper, $f$ is a cofibration if and only if $Uf$ is a 
cofibration
also, i.e. cofibrations are also created by $U$: 
the case of simplicial sets is tautological, and in the case of
simplicial vector spaces that category is an abelian category,
equivalent to the category of $\geq 0$-graded
chain complexes by functors which preserve injections and surjections,
and hence every injection or surjection which induces an equivalence splits.
Additionally, all simplicial surjections are fibrations, but cofibrations are 
characterized as maps having the left lifting property with respect to
fibration equivalences, so cofibrations are simplicial injections also.

\vspace{3mm}
There is also a canonical Quillen
model structure on the category of $\mc$-algebras (with equivalences and
fibrations same as in $s\mb$, i.e. created
by the forgetful functor), and for a cofibrant $\mc$-algebra
$R$, there is a canonical Quillen model structure on the category of
$(\mc,R)$-modules (again with equivalences and
fibrations same as in $s\mb$). This is proven by a ``small object argument",
and the proof works for categories of algebraic structures in $s\mb$ of
very general kinds. (The small object argument is described
in \cite{ds}, 7.12).
We shall refer to representatives of an $s\mb$-equivalence class
as models.


\vspace{3mm}
With the Quillen model structure established, now recall that an
object $X$ is called cofibrant (resp. fibrant) if the map from the initial
object to $X$ (resp. from $X$ to the terminal object) is a cofibration
(resp. fibration). A
{\em cofibrant replacement} (resp. fibrant replacement) of an object $M$ (resp. $N$)
is a map
$$M^{\prime}\r M$$
(resp. $N\r N^{\prime}$) which is a fibration equivalence (resp. cofibration equivalence)
and $M^{\prime}$ (resp. $N^{\prime}$) is cofibrant (resp. fibrant).
We define 
\begin{equation}
\label{e3}
RHom_{(\mc,R)}(M,N)
\end{equation}
as $Hom_{(\mc,R)}(M^{\prime},N^{\prime})$ where $M^{\prime}$ is a cofibrant
replacement of $M$ and $N^{\prime}$ is a fibrant replacement of $N$. 
In our cases, cofibrant and fibrant replacement can be made functorial, 
so \rref{e3} is well defined. Nevertheless, cofibrant and fibrant replacements
are not {\em canonical}, and hence it is appropriate to address the question
of comparing the different $RHom$'s when different selections are made.

\vspace{3mm}
To this end, one uses the following technique. Let $I$ be the standard simplicial
model
of the unit interval. Then we have objects of the form
$$I\otimes M$$
which are cylindrical objects (see \cite{ds}, 4.1)
in the sense that the two maps
$$i_0,i_1:M\r I\otimes M$$
induced by the inclusions of the endpoints to $I$ are equivalences
and $i_{0}\coprod i_{1}:M\coprod M\r I\otimes M$ is a cofibration.
We define a {\em homotopy} of two maps $f_0,f_1:M\r N$ to be a map
$I\otimes M\r N$ which, when composed with $i_j$, gives $f_{j}$.
This is a particular example of what is known as a Quillen left homotopy,
but the present notion has more features which we will find useful. In particular,
the functor $I\otimes ?$ has a right adjoint which we will denote by $F(I,?)$,
and also the internal $Hom$-functor $Hom_{(\mc,R)}$ obeys the relation
\begin{equation}
\label{eobey}
Hom_{(\mc,R)}(I\otimes M,N)\cong Hom_{(\mc,R)}(M,F(I,N))\cong
F(I, Hom_{(\mc,R)}(M,N)).
\end{equation}
To be precise, the $F(I,?)$ on the right hand side of \rref{eobey} is in
$s\mb$ rather than the category of $(\mc,R)$-modules. Note, however, that
since $F(I,?)$ is a limit, the forgetful functor from $(\mc,R)$ to $s\mb$
preserves $F(I,?)$, so it is given by the same construction in $s\mb$ as
in $(\mc,R)$-modules.
All this is formal. Additionally, it is true in our case that $F(I,N)$ is a 
co-cylindrical object (satisfying properties dual to cylindrical object;
also known as path object, see \cite{ds}, 4.12) if
$N$ is fibrant (since, again, this is true in $s\mb$). 
It follows that for any other cofibrant replacement $M^{\prime
\prime}\r M$ there is a comparison map
\begin{equation}
\label{eee1}M^{\prime}\r M^{\prime\prime}
\end{equation}
commuting with the specified maps into $M$, and moreover unique up to homotopy
(in our sense). Hence, by the same principle, we also abtain a map
$$M^{\prime\prime}\r M^{\prime},$$
and the compositions are homotopic to the identity. We call this a homotopy
equivalence of $(\mc,R)$-modules. But then applying $Hom_{(\mc,R)}(?,N^{\prime})$,
we obtain a homotopy equivalence in $s\mb$, which is an equivalence. The treatment
of fibrant replacements is adjoint.

\vspace{3mm}
We shall now turn to the reduction of Theorem \ref{t1} to
another statement.

\vspace{3mm}
\noindent
{\bf Definition:}
Let $S$ be a simplicial set, i.e. a functor $\Delta^{Op}\r Sets$.
Then $S$ can also be viewed as a category $\overline{S}$ with objects
$\cform{\coprod}{n}{}S_{n}$, and morphisms $\phi:s\r t$ where $\phi
\in Mor(\Delta^{Op})$, $s\in S_{n}$ for some $n$, and $\phi(s)=t$.
Let $\Gamma$ be any category. Then a {\em $s\Gamma$-object fibered over $S$} is,
by definition, a functor
$$F:\overline{S}\r \Gamma.$$
We shall write $F_{s}=F(s)$ for $s\in S_{n}$. For a map $i:S\r T$ of simplicial
sets, we have a functor $i^*$ from $s\Gamma$-objects over $T$ to $s\Gamma$-objects
over $S$, given by
$$i^*(F)=F\circ \overline{i}.$$
We shall also make use of the left adjoint to $i^*$, which we shall denote
by $i_{\sharp}$.

Note that $s\Gamma$-objects fibered over the constant simplicial set $*$ are
precisely $s\Gamma$-objects, which helps justify the terminology.

\vspace{3mm}
Specifically, we will now
be interested in the case $\Gamma=s\mb$ as above. 
Clearly, for every pair of objects
$X,Y$ of $s\mb$ fibered over simplicial sets $S,T$, there
is a canonical object $X\boxtimes Y$ fibered over 
$S\times T$. 

\vspace{3mm}
\noindent
{\bf Remark:}
Simplicial realization $|\;|:ss\mb\r s\mb$ is defined as the diagonal functor
\cite{km}:
$$|X|_{n}=X_{nn}.$$
By definition, we have
$$|X\boxtimes Y|=|X|\boxtimes |Y|.$$
We will sometimes drop the symbol $|\;|$ from our notation.


\vspace{3mm}
Let $\mc$ be an operad in simplicial
sets. Then an {\em operad $Q$ in $s\mb$ fibered over $\mc$}
consists of objects $Q(n)$ of $s\mb$ fibered over $\mc(n)$, with 
$\Sigma_n$-action, and unity for $n=1$, and, for each of the compositions
$$\gamma:\mc(k)\times\mc(n_1)\times...\times\mc(n_k)\r
\mc(n_1+...+n_k),$$
a composition
$$\gamma_{\sharp}(Q(k)\boxtimes Q(n_{1})\boxtimes...\boxtimes Q(n_{k}))
\r Q(n_{1}+...+n_{k})$$
satisfying the obvious axioms analogous to the operad axioms \cite{271}.


Now simplicial sets fibered over a simplicial set $S$ are precisely
simplicial sets $X$ over $S$, i.e. arrows $X\r S$. A morphism
in this category is a fibration, cofibration of equivalence if and
only if it has the corresponding property in simplicial sets.
If $U_S$ is the forgetful functor from objects and
morphisms of $s\mb$ fibered
over $S$ to simplicial sets fibered over $S$, then we say that
a morphism $f$ in $s\mb$ fibered over $S$ is a fibration or equivalence
if and only if $U_{S}f$ is a fibration or equivalence. By \cite{ds},
again, this defines
a Quillen model structure on the category of objects of $s\mb$
fibered over $S$.

Finally, on operads (similarly as any type of algebraic
structure) in objects and morphisms of
$s\mb$ fibered over $\mc$ we consider the closed
model structure taking as fibrations (resp. equivalences) sequences
of maps $(A(n)\r\mc(n))_{n}$ which are fibrations (equivalences) in
the category of objects and morphisms of
$s\mb$ fibered over $\mc(n)$.

\vspace{3mm}
\noindent
{\bf Definition:}
Let, for any simplicial set $S$, $h=h_{S}:S\r *$ be the collapse map,
and let 
$$Z=B(\mc(\ell),\mc(1)^{\times \ell},\mc(0)^{\times \ell}).$$
We shall call an operad $Q$ fibered over $\mc$ {\em special}
if, for every $\ell$, the map
\begin{equation}
\label{esp}
B(Q(\ell),Q(1)^{\boxtimes \ell},Q(0)^{\boxtimes \ell})\r h_{Z}^{*}Q(0)
\end{equation}
(where $B$ is the two-sided bar construction),
induced from the composition map
$$h_{Z\sharp}B(Q(\ell),Q(1)^{\boxtimes \ell},Q(0)^{\boxtimes \ell})\r
Q(0)$$
is an equivalence.  

\vspace{3mm}
\noindent
{\bf Remark:}
If $Q$ is fibrant in the category of $s\mb$-operads fibered over $\mc$
such that $\mc(0)=*$, $\mc(1)\simeq *$,
and the unit inclusion
$$i:\{*\}\r\mc(1)$$
is a cofibration equivalence (which we are assuming),
then the counit map
\begin{equation}
\label{emon}
i^{*}Q(1)=h_{\sharp}i_{\sharp}i^{*}Q(1)\r h_{\sharp} Q(1)
\end{equation}
is an equivalence. 
We note that \rref{emon} is a map of monoids in $s\mb$. We shall denote
$$Q_{1}=i^{*}Q(1).$$
Thus, for fibrant operads $Q$ over $\mc$, we can replace $Q(1)$ by $Q_{1}$
in \rref{esp}. Of course, every operad in $s\mb$ over
$\mc$ can be replaced by a fibrant model.
We shall make use of this below.

\vspace{3mm}
Now for a monoid $R$ in $s\mb$, a {\em module} over $R$ is an object
$M$ of $s\mb$ with a map 
$$R\boxtimes M\r M$$
satisfying the usual axioms. Clearly, $R$-modules are precisely algebras
over a monad $C_R$ of the form
$$C_R(X) =R\boxtimes X.$$
Thus, we have a canonical closed model structure on $R$-modules for any
monoid $R$. As before, we define
$$Hom_R(M,N)$$
as the equalizer of the two obvious maps
$$Hom(M,N)\begin{array}{c}\r\\ \r\end{array}Hom(C_{R}M,N),$$
and define $RHom_R(M,N)=Hom_R(M^{\prime},N^{\prime})$
where $M^{\prime}$ is a cofibrant replacement of $M$ and
$N^{\prime}$ is a fibrant replacement of $N$ (with derived independence
on the choice of $M^{\prime}$ and $N^{\prime}$, for $s\mb$
as above).
Note that an example of a $Q_{1}$-module in the preceding
remark is $Q(0)$. 

\vspace{3mm}
Now for two operads $\mc$,$\mathcal{D}$ in simplicial sets, following \cite{BV},
define an operad $\mc\Box\mathcal{D}$ as the quotient of the free operad $\mathcal{F}$
on $\mc\coprod\mathcal{D}=(\mc(n)\coprod\mathcal{D}(n))_{n}$ modulo
identifying the $\mathcal{F}$-operad operations on objects of $\mc$, $\mathcal{D}$ with
the corresponding operations in $\mc$, $\mathcal{D}$, (this includes
units), and the following key relation: for
$\alpha\in\mc(m)$, $\beta\in\mathcal{D}(n)$, 
\begin{equation}
\label{erelation}\alpha(\underbrace{\beta,...\beta}_{\text{$m$ times}})=
\beta(\underbrace{\alpha,...\alpha}_{\text{$n$ times}})\sigma
\end{equation}
where $\sigma$ is a certain permutation reordering terms. To describe this permutation,
consider the ``row by row'' lexicographical bijection
$$\rho_1:\{1,...,m\}\times\{1,...,n\}\r\{1,...,mn\}$$
(i.e. $(11\mapsto 1,12\mapsto 2,...,1n\mapsto n,21\mapsto n+1,...)$), 
and the ``column by column'' lexicographical bijection
$$\rho_2:\{1,...,m\}\times\{1,...,n\}\r\{1,...,mn\}$$
(i.e. $(11\mapsto 1, 21\mapsto 2,..., m1\mapsto m, 12\mapsto m+1,...)$).
The point is that on the left hand side of \rref{erelation}, the entries are
ordered ``row by row'' whereas on the right hand side they are ordered ``column
by column''. Since permutations on operads act on the right, we conclude that
$$\sigma=\rho_2\circ(\rho_1)^{-1}.$$

\begin{proposition}
\label{p1}
There is a canonical map of operads
\begin{equation}
\label{ecube1}
\phi:\mc_k\Box \mc_{\ell}\r \mc_{k+\ell},
\end{equation}
which, in the $n$-th term
$$\phi:(\mc_k\Box\mc_\ell)(n)\r \mc_{k+\ell}(n),$$
is a $\Sigma_n$-equivariant homotopy equivalence
($\Sigma_n$ denotes the symmetric group on $n$ elements). 
An analogous result also holds if
we work in the category of topological spaces (rather than simplicial sets).
\end{proposition}

\vspace{3mm}
\noindent
{\bf Remark:} Note that the Proposition implies that both in the topological
and simplicial cases, the map $\phi$ induces an equivalence of monads
in the sense of \cite{271}.

\vspace{3mm}
We shall prove this at the end of this Section. The second statement
of the Proposition is a theorem of Gerald Dunn (\cite{dunn}, Theorem 2.9).
We give a proof of the first statement, which we need for technical
reasons, at the end of this section. To this end, to fix notation, we
will need to briefly recall the main ideas of Dunn's argument also.
The simplicial extension of Dunn's result contains no substantially new
idea, but due to subtleties must be handled with care.
We will now restate Theorem \ref{t1} as follows:

\begin{theorem}
\label{t2}
Let $Q$ be a special cofibrant fibrant operad in $s\mb$ fibered over $\mc$
where $s\mb$ is as above and $\mc$ is an operad in simplicial sets with
$\mc(0)=*$, $\mc(1)\simeq *$, and such that the $\Sigma_n$-action
on $(\mc\Box\mc_1)(n)$ is free. Then there is
a model of
\begin{equation}
\label{epr*}RHom_{Q_{1}}(Q(0),Q(0))
\end{equation}
which has a natural structure of a $\mc\Box \mc_{1}$-algebra.
\end{theorem}

A discussion is needed to see how Theorem \ref{t2} implies Theorem
\ref{t1}. Recall that in general for a monad $C$ a $C$-functor $D$
(i.e. a functor with a structure map $DC\r D$ satisfying the usual axioms),
and a $C$-algebra $X$, we can define
$$D\otimes_{C}X$$
as the coequalizer of the two maps
$$DCX \begin{array}{c}\r\\ \r\end{array}DX$$
given by the structure maps $DC\r D$, $CX\r X$.
Now for an operad $\mc$ in the category of simplicial sets, 
$\mc$-algebras are algebras of the monad
$$CX=\cform{\coprod}{n\geq 0}{}\mc(n)\otimes_{\Sigma_{n}}X^{\boxtimes n}.$$
We shall define $C$-functors $D_{\ell}$ as follows:
\begin{equation}
\label{E+}
D_{\ell}
X=\cform{\coprod}{n\geq 0}{}\mc(n+\ell)\otimes_{\Sigma_{n}}X^{\boxtimes n}.
\end{equation}
The following statement is obvious upon a moment's reflection:

\begin{proposition}
\label{pop}
Let $\mc$ be an operad in simplicial sets and assume in addition $\mc(0)=*$
(the one point simplicial set). Let
$R$ be a $\mc$-algebra. Then the object
$$A(n)=D_{n}\otimes_{C}R$$
has the natural structure of an operad in $s\mb$ fibered over $\mc$.
\end{proposition}

\Proof
Recall that because the functor $?\boxtimes X$ in $\mb$ has a right adjoint,
$\boxtimes$ is distributive under $\coprod$. Now we need to construct structure
maps 
\begin{equation}
\label{eaux1}
{\diagram
\mc(n+\ell)\otimes_{\Sigma_{n}}R^{\boxtimes n}\boxtimes\\
(\mc(n_{1}+m_{1})\otimes_{\Sigma_{n_{1}}}R^{\boxtimes n_{1}})\boxtimes...
\boxtimes 
(\mc(n_{\ell}+m_{\ell})\otimes_{\Sigma_{n_{\ell}}}R^{\boxtimes n_{\ell}})\dto\\
(\mc(n+n_{1}+...+n_{\ell}+m_{1}+...+m_{\ell})
\otimes_{\Sigma_{n_{1}+...+n_{\ell}}}R^{\boxtimes n+n_{1}+...+n_{\ell}})
\enddiagram}
\end{equation}
compatible with $\otimes_{C}$.
However, realizing that $\otimes$ is nothing but coproducts over various sets,
distributivity applies to this case also, and we see that the map \rref{eaux1}
is obtained just by applying operad composition in $\mc$, and grouping the
$\boxtimes$-powers of $R$. Compatibility with $\otimes_{C}$ is obvious from
operad axioms. Additionally, the fibered structure is obtained
by taking, for an element
$$x\in \mc(n)_{i}$$
($i$ denotes the simplicial index) $A(n)_{x}$ to be the coproduct of all
$$\{y\}\otimes_{\Sigma_{\ell}}R^{\boxtimes \ell}$$
where $y\in \mc(n+\ell)$ is such that 
\begin{equation}
\label{eaux2}
\gamma(y,1,...,1,*,...,*)=x
\end{equation}
where $\gamma$ is operad multiplication, and in \rref{eaux2},
there are $n$ $1$'s and $\ell$ $*$'s.
\qed

The reduction from Theorem \ref{t1} to Theorem \ref{t2} then follows
from the following two results:

\begin{lemma}
\label{pred}
The category of $(\mc,R)$-modules is equivalent to the category of
$h_{\sharp}A(1)$-modules. This equivalence of categories carries
$R$ to $A(0)$. Moreover, if $R$ is a cofibrant $\mc$-algebra
(with $s\mb$ as above), this is a Quillen equivalence, i.e. 
passes on to an equivalence of Quillen homotopy categories.
\end{lemma}

\Proof
The equivalence of categories is established by the fact that both categories
consist of algebras over the same monad:
\begin{equation}
\label{emonad}
X\mapsto h_{\sharp}A(1)\boxtimes X.
\end{equation}
Indeed, it suffices to consider the case when $R$ is a free $\mc$-algebra,
i.e. $R=CX$. Then
$$\begin{array}{l}A(1)=D_{1}\otimes_{C}R=D_{1}X
\\=\cform{\coprod}{n\geq 0}{}
\mc(n+1)\otimes_{\Sigma_{n}}X^{\boxtimes n}.\end{array}$$
Now the free $(\mc,R)$-module on $M$ is
$$\cform{\coprod}{n\geq 0}{}\mc(n+1)\otimes_{\Sigma_{n}}
(X^{\boxtimes n}\boxtimes M)=A(1)\boxtimes M.$$
(By abuse of notation, we treat $h_{\sharp}$ as the forgetful functor,
so this is the same as \rref{emonad}.)

Now if $R$ is a cofibrant $\mc$-algebra, both Quillen model structures on
the respective categories are defined in the same way.
\qed

\vspace{3mm}
\noindent
{\bf Remark:}
Note that this is analogous to a method used by Zhu \cite{zhu}
for vertex operator algebras (see also \cite{dlm}).
Concretely, vertex operator algebras are close to the notion of
algebra over a certain modification to the little $2$-disk operad
$\mathcal{D}$ (see \cite{huang}).
One could elaborate a lot more on that, but 
in the rough analogy, the notion of module over a vertex operator
algebra $V$ corresponds to the notion of $(\mathcal{D},V)$-module.
Now Zhu \cite{zhu} describes an associative algebra $A$ with the
property that, for a rational vertex operator algebra
$V$, irreducible $V$-modules $M$ are in bijective correspondence
with irreducible $A$-modules. However, the algebra $A$ is {\em not}
a precise analogue of $h_{\sharp}A(1)$: in \cite{zhu}, the algebra
$A$ is finite-dimensional, and only acts on the top weight part of $M$. 

\begin{theorem}
\label{tspec}
Let $\mc=\mc_{k}$ be the little cube operad, and let $R$ be a free
$\mc_{k}$-algebra. Let $A$ be defined as above in Proposition \ref{pop}.
Then $A$ is special (although not fibrant).
\end{theorem}

We shall prove this theorem below in Section \ref{sspec}.

\vspace{3mm}
\noindent
{\bf Remark:}
There are other examples of special operads. For example, working in 
spaces (one can get to simplicial sets by applying the singular set functor),
let, for $e\in\mc_{k}(n)$, and a based CW complex $X$,
$$(\Phi_{k}(X))_{e}=Map((e,\partial e), (X,*)).$$
Then there is a standard way to put a topology on
$$\Phi_{k}(X)(n)=\cform{\bigcup}{e\in\mc_{k}(n)}{}\Phi_{k}(X)_{e},$$
making $\Phi_{k}(X)$ an operad fibered over the topological version of
$\mc_{k}$. It can be shown that $\Phi_{k}(X)$ is special if $X$ is
$(k-1)$-connected.

On the other hand, it is easy to construct operads $\mc$ in simplicial
sets such that, for $R=CX$ the associated operad $A$ fibered over $\mc$
is not special: It suffices to take a free operad on a set
(in the category of operads $\mc$ with unit and $\mc(0)=*$).

We shall conclude this section with a 

\vspace{3mm}
\noindent
{\bf Proof of Proposition \ref{p1}:}
We shall first prove the statement for the category of
topological spaces, reviewing essentially
the ideas of Dunn \cite{dunn}. In this case, let $\mc_{k}$ denote the original
topological space models of the little cube operads
rather than the singular set model. In this setting, the map
\rref{ecube1} is obtained by sending a configuration of $n$ little cubes
$$(c_{1},...,c_{n})\in \mc_{k}(n)$$
to
$$(c_{1}\times I^{\ell},...c_{n}\times I^{\ell})\in \mc_{k+\ell}(n),$$
and a configuration of $n$ little cubes
$$(d_{1},...,d_{n})\in\mc_{\ell}(n)$$
to 
$$(I^{k}\times d_{1},...I^{k}\times d_{n}).$$
We will not construct a homotopy inverse of \rref{ecube1} on the whole
$\mc_{k+\ell}(n)$, but instead on a certain subspace 
$\mc_{k+\ell}(n)^{\prime}$ which is weakly equivalent. To define this
subspace, put $m=k+\ell$. We shall call an $n$-tuple of little cubes
\begin{equation}
\label{eaux3a}
e=(e_{1},...,e_{n})
\end{equation}
which is an element of $\mc_{m}(n)$ {\em small} is the following condition
is satisfied: there exists a $p$-tuple of little cubes 
$$f=(f_{1},...,f_{p})$$
forming an element of $\mc_{k}(p)$ and a $q$-tuple of little cubes
$$g=(g_{1},...,g_{q})$$
forming an element of $\mc_{\ell}(q)$ such that every little cube $e_{i}$
lies in the interior of precisely one set
\begin{equation}
\label{eaux3}
f_{h}\times g_{j},
\end{equation}
and every set \rref{eaux3} contains at most one little cube $e_{i}$.
The space $\mc_{m}(n)^{\prime}$ is the subspace of $\mc_{m}(n)$ consisting
of precisely all small $n$-tuples. Now we claim that the inclusion
$\mc_{m}(n)^{\prime}\subset \mc_{m}(n)$ is a homotopy equivalence.
Indeed, let $e$ be as in \rref{eaux3a}. Then, for $\lambda\in (0,1]$,
define $\lambda e$ as the little cube configuration obtained by scaling
each little cube $e_{i}$ by a factor $\lambda$ in its center. Then 
we know that $\lambda e\in \mc_{m}(n)$, and if $e$
is an element of $\mc_{m}(n)^{\prime}$, then so is $\lambda e$.
Furthermore, it is easy to see that for every $e$ there exists a $\lambda$
such that $\lambda e\in\mc_{m}(n)^{\prime}$ (the statement is true trivially
if every little cube is replaced by one point, namely its center). By the
Lebesgue number theorem, for every compact subset 
$K\subset \mc_{m}(n)$,
there exists a $\lambda\in (0,1]$ such that 
$\lambda K\subset \mc_{m}(n)^{\prime}$.
Furthermore, $t.Id$, $t\in[\lambda,1]$ is a homotopy between 
$K$ and $\lambda K$,
which moreover stays in $\mc_{m}(n)^{\prime}$ if 
$K\subset \mc_{m}(n)^{\prime}$.
Clearly, this construction is $\Sigma_n$-equivariant.
By the Whitehead theorem, the inclusion
$\mc_{m}(n)^{\prime}\subset \mc_{m}(n)$ is a $\Sigma_n$-weak 
equivalence, and
hence in fact a $\Sigma_n$-homotopy equivalence, 
since $\mc_{m}(n)$ is a $\Sigma_n$-CW complex.

\vspace{3mm}
Now we will construct a right inverse $\psi$ to the map $\phi$ when
restricted to $\mc_{m}(n)^{\prime}$. In fact, this construction is obvious:
simply compose $f$ with $p$ copies of $g$, and each entry with either $*$ or
the appropriate elements of $\mc_{k}(1)$, $\mc_{\ell}(1)$ to ensure that
$\phi\psi=Id_{\mc_{m}(n)^{\prime}}$. It should be noted that the map
is well defined and continuos, because its value on an element $e$ does
not depend on the choice of $f$, $g$: any two 
choices have a ``common subdivision'',
which produces the same element by the fundamental relation
\rref{erelation}. To be more precise, if $A$, $B$ are two sets of disjoint
little cubes in $I^{k}$, the {\em common subdivision} of $A$, $B$ is
$$\{a\cap b|a\in A, b\in B\}.$$
(While ordering of the cubes of course matters in the operad structure,
we do not have to specify it in this definition, as any two ordering are
related by the symmetric group action, and hence any ordering will do.)

\vspace{3mm}
Note that we are not yet done: We must still produce a homotopy left inverse
to $\phi$. But now let 
$$(\mc_{k}\Box \mc_{\ell})(n)^{\prime}=\phi^{-1}(\mc_{m}(n)^{\prime}).$$
First of all, note that the inclusion 
$$(\mc_{k}\Box \mc_{\ell})(n)^{\prime}\subset (\mc_{k}\Box \mc_{\ell})(n)$$
is a weak equivalence by the same argument as above: we may emulate
the homotopy corresponding to multiplying $e$ by $t$ by
composing an element of $(\mc_{k}\Box \mc_{\ell})(n)$ with $n$ copies of
$\gamma(t.1_{k},t.1_{\ell})$ where $\gamma$
is operad composition, and $1_{k}\in \mc_{k}(1)$, $1_{\ell}\in
\mc_{\ell}(1)$ are the unit elements. So we are done if we can show that
$\psi\phi=Id$ on $(\mc_{k}\Box \mc_{\ell})(n)^{\prime}$.
But this is just a refinement of the above argument that $\psi$ did not
depend on the choice of $f$, $g$: one may form common subdivisions
with the $\mc_{k}$ and $\mc_{\ell}$ elements $u_{i}$ figuring in the definition
of an element of $(\mc_{k}\Box \mc_{\ell})(n)^{\prime}$, and use the
relation \rref{erelation} to show that the common subdivision
produces the same element as using either $f$, $g$, or $u_{i}$.
Again, everything is $\Sigma_n$-equivariant.
This concludes the proof of our statement for the category of topological
spaces.

\vspace{3mm}
We shall now study what changes when we work in the category of simplicial
sets, using the singular
sets of $\mc_{k}(n)$ etc. instead of the actual spaces. In
this part of the proof,
we will find it convenient to display the 
simplicial set functor $S$ explicitly,
to prevent confusion. 
Much of the idea is the same. For example, the construction of the map $\phi$
is got simply by applying the singular set functor to the space level $\phi$
(it is useful to note that there is always a canonical
map $S\mc\Box S\mathcal{D}\r S(\mc\Box \mathcal{D})$). However, when 
constructing the map $\psi$, we must adapt the definition of
$\mc_{m}(n)^{\prime}$. In fact, we must introduce the notion of {\em small
singular simplex} in $S\mc_{m}(n)$ as follows: if we represent the singular
simplex by an $n$-tuple of singular simplices
$$(e_{1},...,e_{n})$$
in the space of little cube, (i.e. for $t\in\Delta_{N}$, for some $N$,
$e_{i}(t)$ is a little cube), then there must exist a {\em uniform} (i.e.
independent of $t$) choice of $f$, $g$ such that $e_{i}(t)$ satisfy the
above condition in place of $e_{i}$ for each $t$. This means, roughly, that the value
of a small singular simplex at each $t$ is required to be small, but 
also the values of the singular simplex must vary only by a ``small'' amount.
We take $S\mc_{m}(n)^{\prime}$ as the simplicial set of small singular simplices
in $S\mc_{m}(n)$. Then the map $\psi$ may be define completely analogously as
in the case of spaces, by passing to singular sets.

\vspace{3mm}
It is, further, correct to think of $\psi$ as a 
right homotopy inverse to $\phi$,
as it can be shown that $S\mc_{m}(n)^{\prime}\subset S\mc_{m}(n)$ is an 
equivalence:
this is a special case of a general theorem stating that for any 
open covering $(U_{i})$
of a space $X$, the inclusion of the
sub-simplicial set of the singular set of $X$ consisting
of singular simplices whose images are in one of the $U_{i}$'s is
an equivalence. (In our case, $X$ is the topological $S\mc_{m}(n)^{\prime}$.)

\vspace{3mm}
However, we must still find a left homotopy inverse to $\phi$. 
By barycentric subdivision of a simplicial
set $T$ we shall mean the simplicial set which is the classifying
space of the category whose objects are non-degenerate
elements of $T$ and morphisms are faces.
Denote the $j$-fold 
barycentric subdivision of a simplicial set
$T$ by $T^{(j)}$. 
It is a well known fact that
if, for the moment, $|?|$ denotes topological simplicial realization,
then there is a canonical homeomorphism
$$|T|\cong |T^{(j)}|.$$
On the other hand,
it is a standard fact (used for example in proving
Eilenberg-Steenrod axioms for singular homology) that we have a canonical
simplicial map
\begin{equation}
\label{eaux10}
\iota_{j}:S^{(j)}X\r SX
\end{equation}
for any space $X$, sending the ``algebraic barycenter'' to the ``topological
barycenter'' (obviously, it suffices to consider $j=1$). 
Now the simplicial operad $S\mc_k\Box S\mc_{\ell}$
is not (at least a priori) the singular set of any space,
although, as pointed above, there is a canonical map
$$S\mc_k\Box S\mc_{\ell}\r S(\mc_k\Box\mc_{\ell}).$$
However, we claim that there nevertheless is a
canonical simplicial map
\begin{equation}
\label{efff}{\tilde{\iota}_j:(S\mc_k\Box S\mc_{\ell})^{(j)}\r S\mc_k\Box S\mc_{\ell}}
\end{equation}
constructed in the same way as the map $\iota_j$
for singular sets: the point is that
any word in singular simplices making up an element of 
$(S\mc_k\Box S\mc_{\ell})_{n}$ 
can be written (and any relation among such elements remains
valid) on restrictions of each of those singular simplices
to the same simplex 
of an (iterated) barycentric
subdivision of the standard simplex. We further see that by the
same argument as for spaces, upon applying simplicial realization, 
\rref{efff} becomes homotopic to the canonical homeomorphism between
the realization of a simplicial set and its iterated barycentric subdivision.
Hence \rref{efff} is an equivalence. 

Thus, we would be done if we could
show that $\phi\circ\tilde{\iota}_j$ applied
to the $n$-th space of the source operads lands in $S\mc_{m}(n)^{\prime}$
for some $j$. While this is obviously too much to expect, it is however
true that for any {\em finite} simplicial subset $T$ 
of $S\mc_{k}\Box S\mc_{\ell}$
there exists a $j$ such that the simplicial subset $T^{(j)}$ of 
$$(S\mc_{k}\Box S\mc_{\ell})^{(j)}$$ 
consisting of $j$-fold barycentric
subdivisions of simplices of $T$ does have the property that
$$\phi\circ\tilde{\iota}_j(T^{(j)})\subset (S\mc_{m})(n)^{\prime}.$$
By the Whitehead theorem, this is enough.
\qed

\section{Lax algebras}
\label{srect}
\label{s3}

This story is complex enough that it seems worth telling for
the category of (simplicial) sets first. We use essentially the ideas
of \cite{271}, \cite{km}, but as far as we know, they have not
been recorded with the specific nuances needed here. 

Consider the category $\Sigma$ whose objects are finite sets
$\{1,...,n\}$, $n\geq 0$, and morphisms are permutations. Then
an operad is, in particular, a functor on $\Sigma$.

\vspace{3mm}
Now for any operad $\mc$, a {\em lax $\mc$-algebra} is a category
$\Gamma$ where each $n$-ary operation $\mu$ of $\mc$ corresponds to a functor
\begin{equation}
\label{ena+}
{\begin{array}{c}
\underline{\mu}:\Gamma^{\times n}\r \Gamma,\\
\end{array}}
\end{equation}
Further, for each $\mu\in \mc(n)$, $\mu_{i}\in\mc(n_{i})$, $i=1,...,n$,
we are required to have natural {\em coherence morphisms}
\begin{equation}
\label{ecoh+}
\phi:\umu(\umu_{1},...,\umu_{n})\r
\underline{\gamma(\mu,\mu_{1}...,\mu_{n})},\;
\phi:\underline{1}\r 1
\end{equation}
where $\gamma$ is the composition in $\mc$, and also
for $\kappa:\{1,...,n\}\r \{1,...,n\}\in Mor(\Sigma)$,
\begin{equation}
\label{ecoh+1}
\phi:\underline{\kappa\mu}\cong \underline{\mu}\kappa^{*}
\end{equation} 
where $\kappa^*:\Gamma^{\times m}\r \Gamma^{\times n}$ 
is the categorically induced
map by $\kappa$.
In addition, the coherence 
isomorphisms \rref{ecoh+}, \rref{ecoh+1} are required to 
satisfy {\em coherence diagrams} expressing a cocycle condition for the 
isomorphisms \rref{ecoh+}. This means that for a $3$-fold composition,
applying the coherence isos in either order gives the same result:
$$
\diagram
&\umu(\umu_{1}(\umu_{11},...,\umu_{1m_{1}}),...,\umu_{n}(\umu_{n1},...,
\umu_{nm_{n}}))\dlto^{\phi}\drto^{\phi}&\\
\protect
{\begin{array}{c}
\umu(\underline{\gamma(\mu_1,\mu_{11},...,\mu_{1m_{1}})},\\...,
\underline{\gamma(\mu_n,\mu_{n1},...,\mu_{nm_{n}})})\end{array}}
\drto^{\phi}
&&
\protect
{\begin{array}{c}
\underline{\gamma(\mu,\mu_{1},...,\mu_{n})}\\(\umu_{1},...,\umu_{nm_{n}})
\end{array}}
\dlto^{\phi}\\
&\underline{\gamma(\mu,\gamma(\mu_{1},\mu_{11},...,\mu_{1m_{1}}),...,
\gamma(\mu_{n},\mu_{n1},
...,\mu_{nm_{n}}))}.&
\enddiagram
$$
The diagrams coming from operad unit axioms are
$$\diagram
\underline{1}(\umu)
\rto^{\phi}\dto_{\phi(\umu)}&\underline{\gamma(1,\mu)}\\
\umu,\urto_{=}&
\enddiagram$$
$$\diagram
\umu(\underline{1},...,\underline{1})\rto^{\phi}\dto_{\umu(\phi)}&
\underline{\gamma(\mu,1,...,1)}\\
\umu.\urto_{=}&
\enddiagram
$$
Regarding equivariance, there are two diagrams (associativity and unit)
coming from the axiom that $\mathcal{C}:\Sigma\r Sets$ is a functor. If
$\kappa:\{1,...,n\}\r\{1,...,n\}$, $\lambda:\{1,...,n\}\r\{1,...,n\}$
and $\iota=Id:\{1,...,n\}\r\{1,...,n\}$, the diagrams are
$$
\diagram
\underline{\lambda\kappa\mu}\rto^{\cong\phi}\dto_{\cong\phi}&\underline{\kappa
\mu}\lambda^{*}\dto^{\cong\phi\lambda^*}\\
\umu(\lambda\kappa)^*\rto^= &\umu\kappa^*\lambda^*,
\enddiagram
$$
$$
\diagram
\underline{\iota\mu}\rto^{\cong\phi}\dto_=&\umu\iota^*\dlto^=\\
\umu.&
\enddiagram
$$
Finally, one diagram comes from the commutation relation between composition
and equivariance in an operad. Let $\lambda:\{1,...,n\}\r\{1,...,n\}$,
$\kappa_i:\{1,...,k_i\}\r\{1,...,k_i\}$, $i=1,...,n$. Then there is
a permutation 
$$\lambda\wr(\kappa_1,...,\kappa_n):\{1,...,\cform{\sum}{i=1}{n}
k_{i}\}
\r\{1,...,\cform{\sum}{i=1}{n}k_{i}\}$$
given by
$$(\lambda\wr(\kappa_1,...,\kappa_n))(\cform{\sum}{i=1}{j-1}k_{\lambda(i)}+p)=
\cform{\sum}{i=1}{\lambda(j)-1}k_{i}+\kappa_{\lambda(j)}(p)$$
for $j=1,...,n$, $p=1,...k_{\lambda(j)}$. The diagram then reads
$$\diagram
\underline{(\lambda\wr(\kappa_1,...,\kappa_m))\gamma(
\mu,\mu_{\lambda(1)},...,\mu_{\lambda(n)})}
\rto^{\phi} & \underline{\gamma(\mu,\mu_{\lambda(1)},...,
\mu_{\lambda(n)})}(\lambda\wr(\kappa_1,...,
\kappa_m))^{*}\\
\underline{\gamma(\lambda\mu,\kappa_{1}\mu_1,...\kappa_{m}\mu_{m})}
\uto_{=}&\umu(\umu_{\lambda(1)},...,\umu_{\lambda(n)}
)(\lambda\wr(\kappa_1,...,\kappa_m))^{*}
\uto^{\phi(\lambda\wr(\kappa_1,...,\kappa_m))^*}\\
\underline{\lambda\mu}(\underline{\kappa_{1}\mu_1},...,\underline{\kappa_{
m}\mu_{m}})\uto_{\phi}\rto_{\phi(\phi,...,\phi)}&
\umu\lambda^*(\umu_1\kappa^{*}_{1},...,\umu_{m}\kappa^{*}_{m}).\uto^{=}
\enddiagram
$$

In another formulation, we may consider the free operad $Op(\mc)$
on $\mc$ (i.e. the operad obtained by iterating the operations in $\mc$,
and performing substitutions and insertions of unit),
and assign a canonical iso to any two elements in $Op(\mc)(n)$
whose images in $\mc(n)$ coincide: these isos are required to be transitive
and compatible with substitution and composition.
Note that these isos, together with their iterations,
now make each $Op(\mc)(n)$ into a category $G(n)$,
whose fiber over each object of $\mc$
has a final object. 
Then $Op(\mc)$ defines a $2$-monad
in the sense of Blackwell, Kelly and Power \cite{bkp}.
Lax $\mc$-algebras can then be identified with $2$-algebras
over the $2$-monad $Op(\mc)$. 

\vspace{3mm}
Now let $\mc$ be an operad. We define a simplicial operad $\overline{\mc}$ 
as
$$\begin{array}{l}\overline{C}(n)=BG(n)=
B_{Obj(G(n))}(Obj(G(n)),Mor(G(n)),Obj(G(n)))
\\=B_{Op(C)(n)}(Op(C)(n),Mor(G(n)),Op(C)(n)).\end{array}$$
We have a canonical projection
$$|\overline{\mc}|\r\mc$$ 
which is an equivalence since the fiber of
$G(n)$ over every element in $\mc$
has a final object.

\begin{proposition}
\label{prect1}
Let $\Gamma$ be a lax $\mc$-algebra. Then $B\Gamma$ is canonically a
$\overline{\mc}$-algebra.
\end{proposition}

\Proof
Let
$$A=\left(
\begin{array}{ccc}
a_{11}&...&a_{1n}\\
...&...&...\\
a_{m1}&...&a_{mn}
\end{array}
\right)
$$
be a matrix of morphisms in $\Gamma$, $Ta_{ij}=Sa_{i,j+1}$. ($S,T$ mean source
and target.)
Let $\phi=(\phi_{1},...,\phi_{n})$ be a composable $n$-tuple in $G(m)$,
$$\phi_{i}:\sigma_{i-1}\r\sigma_{i}.$$
Then
$$\phi A=(
\phi_{1}\sigma_{0}
\left(\begin{array}{c}
a_{11}\\...\\a_{m1}
\end{array}\right)
,
\phi_{2}\sigma_{1}
\left(\begin{array}{c}
a_{12}\\...\\a_{m2}
\end{array}\right)
,
...,
\phi_{n}\sigma_{n-1}
\left(\begin{array}{c}
a_{1n}\\...\\a_{mn}
\end{array}\right)
).
$$
\qed

\vspace{3mm}
In an arbitrary symmetric monoidal ground category $\mb$
with monoidal structure $\boxtimes$, 
we can analogously speak of $\mc$-algebra for an operad $\mc$ (see also
Ginzburg-Kapranov \cite{gk}). 

Specifically, we introduce the
following terminology: a {\em category $C$ enhanced in $\mb$}
consists of objects 
$$Obj(C),Mor(C)\in Obj(\mb),$$ 
and morphisms
$S,T,Id,\gamma\in Mor(\mb)$, $S,T:Mor(C)\r Obj(C)$, $Id:Obj(C)\r Mor(C)$,
$\gamma:Mor(C)\prod_{Obj(C)} Mor(C)\r Mor(C)$ (the source of $\gamma$
denotes a pullback), with the usual axioms. 

\vspace{3mm}
\noindent
{\bf Remark:} Note that this differs from a more usual notion
of {\em $\mb$-enriched} category $\Gamma$ in which case we have
a set of objects, and for objects $X,Y\in Obj(\Gamma)$,
an object $\underline{Mor}_{\Gamma}(X,Y)\in Obj\mb$,
with morphisms
$$\underline{Mor}_{\Gamma}(X,Y)\boxtimes\underline{Mor}_{\Gamma}
(Y,Z)\r \underline{Mor}_{\Gamma}(X,Z).$$
Hence, the different terminology. Enriched categories, however,
will also be relevant in the next section.

\vspace{3mm}
We then have an object
$$BC\in Obj(s\mb)$$
defined by
$$BC_{k}=\underbrace{Mor(C)\prod_{Obj(C)}...\prod_{Obj(C)}
Mor(C)}_{\text{$k$ times}}.$$

Now let $\mc$ be an operad in (simplicial) sets. There is a notion of 
$\mc$-algebra $X$ enhanced over $\mb$. The structure maps are of the form
\begin{equation}
\label{eenr+}
\gamma:\mc(n)\otimes (\underbrace{X\boxtimes...\boxtimes X}_{\text{$n$ times}}
)\r X
\end{equation}
with the usual diagrams mimicking the diagrams defining a $\mc$-algebra.

Now a {\em lax $\mc$-algebra enhanced over $\mb$} is
a category $C$ enhanced over $\mb$ together with functorial structure
maps \rref{eenr+} for $X=Obj(C), Mor(C)$ and coherence 
morphism structure of the following
form:

For $x\in \mc(n)$, $y_{i}\in\mc(k_{i})$, $i=1,...,n$, 
$$\phi_{x,y_{1},...,y_{n}}
\in Mor(\mb):Obj(C)^{\boxtimes \sum k_{i}}\r Mor(C)$$
with a commutative diagram
$$
\diagram
Obj(C)^{\boxtimes n}\dto_{\gamma_{x}}&Obj(C)^{\boxtimes \sum k_{i}}
\lto_{\boxtimes \gamma_{y_{i}}}\dto^{\gamma_{\gamma(x,y_{1},...,y_{n})}}\\
Obj(C)\rto^{\phi_{x,y_{1},...,y_{n}}}& Obj(C)
\enddiagram
$$
and corresponding coherence diagrams, mimicking the coherence diagrams
in the case of sets. One also must not forget to include 
coherence isomorphisms corresponding to
permutations (the $\Sigma_{n}$-action on the
$n$-th space of an operad), and coherence diagrams corresponding
to axioms involving composition and equivariance (cf. \cite{271}).
The fact that ``the targets of $\phi_{x,y_{1},...,
y_{n}}$ are iso" is expressed, for example, by giving an `inverse' map
$$\psi_{x,y_{1},...,y_{n}}\in Mor(\mb):
Obj(C)^{\boxtimes \sum k_{i}}\r Mor(C).$$
Thus, we get

\begin{proposition}
\label{prect1a}
Let $\Gamma$ be a lax $\mc$-algebra enhanced over
$\mb$. Then $B\Gamma$ is canonically a
$\overline{\mc}$-algebra in $\mb$.
\end{proposition}
\qed

\vspace{5mm}
\section{Lax algebras enhanced over categories}
\label{s4}

\vspace{3mm}
We shall now need to consider an even further generalization. Let
us, again, first work in the context of simplicial sets, where the structure
is simpler than in the $\mb$-enhanced case. In this case, the appropriate
notion is a lax $\mc$-algebra $\Gamma$ enhanced over categories. This means
that $\Gamma$ has a structure analogous to that of lax $\mc$-algebra,
where both $Obj(\Gamma)$, $Mor(\Gamma)$ are categories 
(which we will refer to as
the {\em vertical} categories), and all structure maps of lax $\mc$-algebra
are functors (and coherence diagrams commute strictly, rather than just
up to natural isomorphisms). Therefore, in addition to the vertical categories,
we get two lax $\mc$-algebras
\begin{equation}
\label{ellax1}
(Obj_{Obj(\Gamma)}, Obj_{Mor(\Gamma)})
\end{equation}
and
\begin{equation}
\label{ellax2}
(Mor_{Obj(\Gamma)}, Mor_{Mor(\Gamma)}).
\end{equation}
(We use the subscript notation to distinguish this structure from the
vertical categories.) In particular, then, \rref{ellax1}, \rref{ellax2}
are categories, and we will refer to them as the {\em horizontal} categories.
To spell out our notation completely, the vertical categories then are
$$ Obj(\Gamma)=(Obj_{Obj(\Gamma)}, Mor_{Obj(\Gamma)}),$$
$$ Mor(\Gamma)=(Obj_{Mor(\Gamma)}, Mor_{Mor(\Gamma)}).$$
Now suppose we are given a lax $\mc$-algebra $\Gamma$ enhanced over categories.
Then performing the horizontal bar construction (which we will denote by $B_{hor}$),
we obtain a strict simplicial $\overline{\mc}$-algebra enhanced over 
categories (using the vertical structure):
$$B_{hor}(\Gamma)=\left(\begin{array}{l}B_{hor}(Obj_{\Gamma})\\
B_{hor}(Mor_{\Gamma})\end{array}\right).$$
Observe that $B_{hor}(\Gamma)$ is ``almost'' a $\mc\Box\mc_{1}$-algebra:
the $\mc$-algebra structure was just described, and the ``$\mc_{1}$-structure''
can be pulled back from the (associative) categorical composition. The 
difficulty with that is that the source of the categorical composition
is not the product
\begin{equation}
\label{ellax3}
B_{hor}(Mor_{\Gamma})\times B_{hor}(Mor_{\Gamma}),
\end{equation}
but the fibered product
\begin{equation}
\label{ellax4}
B_{hor}(Mor_{\Gamma})\times_{B_{hor}(Obj_{\Gamma})} B_{hor}(Mor_{\Gamma}).
\end{equation}
There is a natural inclusion of \rref{ellax4} to \rref{ellax3}.
We need a technique for extending the domain of the composition product
from \rref{ellax4} to \rref{ellax3}.

The technique we shall use is the two-sided bar construction of
monads \cite{271}. The ground category is the category $\mathcal{G}$
of {\em graphs} of $\overline{\mc}$-algebras over $B=B_{hor}(Obj_{\Gamma})$, 
which means $\overline{\mc}$-algebras $X$ with two maps of
$\overline{\mc}$-algebras 
$$S,T:X\r B.$$
Then the monad in $\mathcal{G}$ which defines categories in $\omc$
with objects $B$ (where, as above, composition commutes with 
$\omc$-algebra structure) is
\begin{equation}
\label{ellax5}
DX=\cform{\coprod}{n\geq 0}{}\underbrace{X\times_{B}...
\times_{B}X}_{\text{$n$ times}}.
\end{equation}
On the other hand, the monad which defines monoids in $\omc$-algebras
(i.e. again, we require that the $\omc$-algebra structure commutes with
composition) is
\begin{equation}
\label{ellax6}
EX=\cform{\coprod}{n\geq 0}{}\underbrace{X\times...
\times X}_{\text{$n$ times}}.
\end{equation}
Clearly, we have a map of monads $D\r E$, and we can therefore consider the
$2$-sided bar construction of monads
\begin{equation}
\label{ellax7}
B(E,D,X).
\end{equation}
Then (the realization of) \rref{ellax7} is a monoid in $\omc$-algebras,
i.e. a $\omc\Box\mc_{1}$-algebra by pullback. Of course, we would like
to compare the homotopy type of \rref{ellax7} to the homotopy type of
$X$. As usual, we have the comparison map
\begin{equation}
\label{ellax7a}
{\diagram
X& B(D,D,X)\lto_{\simeq}\rto & B(E,D,X).
\enddiagram}
\end{equation}
When is the second map \rref{ellax7a} an equivalence? An obvious condition is
\begin{equation}
\label{ellax8}
B=B_{hor}(Obj_{\Gamma}) \;\text{is contractible.}
\end{equation}
However, \rref{ellax8} per se unfortunately does not suffice, as we need
some local condition. In order to formulate the condition,
we must slightly change our context: we shall 
actually assume that \rref{ellax1},
\rref{ellax2} are {\em simplicial} categories where the simplicial
structure of \rref{ellax1} is constant (we will see in the next section that
such situation arises naturally). We may of course 
always realize to make objects
simplicial, but for the purposes of the following condition, $X$, $B$ are
then bisimplicial sets, where in one of the simplicial directions $B$ 
is constant.
The condition reads as follows:
\begin{equation}
\label{ellax8a}
\parbox{3.5in}{Considering the projection map $S\times T:X\r B\times B$,
then for every $b\in (B\times B)_{n}$, and every face $d_{i}$ (resp.
$s_{i}$, the
face (degeneracy) $(S\times T)^{-1}b\r (S\times T)^{-1}d_{i}b$,
$(S\times T)^{-1}b\r (S\times T)^{-1}s_{i}b$ are equivalences.}
\end{equation}

\begin{lemma}
\label{interm}
Assuming \rref{ellax8a}, every projection map
\begin{equation}
\label{ellax9}
\underbrace{X\times_{B}...\times_{B}X}_{\text{$n$ times}}\r B^{\times (n+1)}
\end{equation}
is a quasifibration.
\end{lemma}

\Proof
It suffices to prove that $S\times T:X\r B\times B$
is a quasifibration. The Bousfield-Friedlander theorem (Theorem B.4
in \cite{BF}) implies that such map is a quasifibration provided
that, in their terminology,

(i) the vertical homotopy groups $\pi_{n}^{v}(X,x)$ satisfy
the $\pi_*$-Kan condition for all $n\geq 1$

(ii) the simplicial map $p:\pi_{0}^{v}(X)\r B\times B$
is a Kan fibration. 

But these assumptions follow from our condition \rref{ellax8a}.
\qed

So, if $F$ is any fiber of 
\rref{ellax9}, we have a diagram
$$\diagram
&X\times_{B}...\times_{B}X\ddto\\
F\urto^{\simeq}\drto_{\simeq}&\\
&X\times...\times X
\enddiagram
$$
and hence the vertical map is an equivalence. It follows that the right hand
map \rref{ellax7a} is an equivalence. Note also that this implies that $X$ is
equivalent, via inclusion, to any fiber of any of the maps $S,T$.

Now we need to translate \rref{ellax8a} to some condition on categories
which could be applied in the case when $X=B_{hor}(Mor_{\Gamma})$. The following
condition is obviously sufficient ($Hom_{vert}$ denotes vertical
$Hom$ sets):
\begin{equation}
\label{ellax8b}
\parbox{3.5in}{Let $M,N,M^{\prime},N^{\prime}\in Obj_{Obj(\Gamma)}$,
let $f:M\r M^{\prime}$, $g:N\r N^{\prime}$ be horizontal morphisms.
Then the natural projections $S:Hom_{vert}(f,g)\r Hom_{vert}(M,N)$,
$T:Hom_{vert}(f,g)\r Hom_{vert}(M^{\prime},N^{\prime})$ are equivalences,
and $S$ is a fibration.}
\end{equation}
We shall call $\Gamma$ {\em distinguished} if \rref{ellax8}, \rref{ellax8b}
are satisfied. We therefore have proven

\begin{proposition}
\label{pllax*}
Let $\Gamma$ be a distinguished lax $\mc$-algebra 
enhanced over categories. Then
for each $M,N\in Obj_{Obj(\Gamma)}$, $Hom_{vert}(M,N)$ are naturally equivalent
to each other,
and moreover naturally equivalent
to a $\omc\Box\mc_{1}$-algebra.
\end{proposition}

\vspace{3mm}
\noindent
{\bf Remark:}
Recall that May's two-sided bar construction of monads \cite{271}
allows us, for any map of operads in simplicial
sets $\mathcal{D}_1\r\mathcal{D}_2$ which is
an equivalence and such that for all $n$ the
action of $\Sigma_n$ on $\mathcal{D}_i(n)$
is free, and for any $\mathcal{D}_1$-algebra $X$, to construct
the equivalent $\mathcal{D}_2$-algebra
$$B(D_2,D_1,X)$$
where $D_i$ are the monads associated with the operads $\mathcal{D}_i$.
Therefore, if in Proposition \ref{pllax*}
$\mc\Box \mc_1(n)$ is $\Sigma_n$-free,
so is $\overline{\mc}\Box \mc_1(n)$
(since it maps into $\mc\Box \mc_1(n)$), and we are allowed to replace 
$\omc\Box\mc_1$ with $\mc\Box\mc_1$.

\vspace{3mm}
Now to treat the case enhanced over a symmetric monoidal category
$\mb$, we start by defining a {\em category $\Gamma$ enhanced
in $\mb$-categories}. Such structure consists of the following data:
First, we have a `horizontal object category' 
$$Obj_{\Gamma}=(Obj_{Obj(\Gamma)}, Obj_{Mor(\Gamma)})$$
which is an ordinary category (i.e. ``enhanced'' only over Sets).

Next, we have a `vertical object category' given by specifying, for each
$x,y\in Obj_{Obj(\Gamma)}$, a
$$Hom(x,y)\in Obj(\mb).$$
Further, there are specified maps
$$Id\in Mor(\mb):1_{\Box }\r Hom(x,x),$$
$$\gamma\in Mor(\mb):Hom(x,y)\boxtimes Hom(y,z)\r Hom(x,z)$$
with usual axioms of associativity and unity and a `vertical morphism category'
specifying, similarly, for 
$f,g\in Obj_{Mor(\Gamma)}$, a
$$Hom(f,g)\in Obj(\mb).$$
Further, there are specified maps
$$Id\in Mor(\mb):1_{\Box }\r Hom(f,f),$$
$$\gamma\in Mor(\mb):Hom(f,g)\boxtimes Hom(g,h)\r Hom(f,h)$$
with usual axioms of associativity and unity. (Note: we see
that the vertical category is in fact enriched over $\mb$
in the usual sense - see the remark in the previous section.)

Next, we have a `horizontal morphism category' 
$Mor_{\Gamma}$ enhanced over $\mb$:
For 
$$\diagram x_{0}\rto^{f_{1}}&x_{1}\rto^{f_{2}}&x_{2}\enddiagram,$$
$$\diagram y_{0}\rto^{g_{1}}&y_{1}\rto^{g_{2}}&y_{2}\enddiagram$$ 
in $Obj_{\Gamma}$,
morphisms $T\in Mor(\mb):Hom(f_{1},g_{1})\r Hom(x_{1},y_{1})$,
$S\in Mor(\mb): Hom(f_{2},g_{2})\r Hom(x_{1},y_{1})$ and
$$\gamma\in Mor(\mb):Hom(f_{1},g_{1})\prod_{Hom(x_{1},y_{1})}
Hom(f_{2},g_{2})\r Hom(f_{2}f_{1},g_{2},g_{1})$$
and also
$$Id\in Mor(\mb):Hom(x,x)\r Hom(Id_{x},Id_{x}).$$
Finally, there is a diagram of commutativity between vertical and
horizontal composition. This diagram expresses the equality of two
maps from
\begin{equation}
\label{esource}
(Hom(f_{1},g_{1})\prod_{Hom(x_{1},y_{1})}Hom(f_{2},g_{2}))
\boxtimes (Hom(g_{1},h_{1})\prod_{Hom(y_{1},z_{1})}Hom(g_{2},h_{2}))
\end{equation}
to
\begin{equation}
\label{etarget}
Hom(f_{2}f_{1},h_{2}h_{1}).
\end{equation}
The first map maps \rref{esource} to
$$Hom(f_{2}f_{1},g_{2}g_{1})\boxtimes Hom(g_{2}g_{1},h_{2}h_{1})$$
by horizontal composition, and then maps to \rref{etarget} by
vertical composition. The second map maps \rref{esource} to
$$(Hom(f_{1},g_{1})\boxtimes Hom(g_{1},h_{1}))
\prod_{(Hom(x_{1},y_{1})\boxtimes Hom(y_{1},z_{1}))}
(Hom(f_{2},g_{2})\boxtimes Hom(g_{2},h_{2}))$$
using the limit properties of a pullback, followed by
a map to
$$Hom(f_{1},h_{1})\prod_{Hom(x_{1},z_{1})} Hom(f_{2},h_{2})$$
by vertical composition, and then to \rref{etarget} by horizontal
composition.

\vspace{3mm}
Now the axioms of a {\em lax $\mc$-algebra $\Gamma$ enhanced in 
$\mb$-categories} consists of the following data:
\vspace{2mm}
\begin{enumerate}
\item
An ordinary lax $\mc$-algebra structure on $Obj_{\Gamma}$.

\item
\label{i2}
A structure of lax $\mc$-algebra enhanced in $\mb$ on $Mor_{\Gamma}$,
compatible with the fibering of $Mor_{\Gamma}$ over $Obj_{\Gamma}$.

\item
Compatibility diagrams of \rref{i2} with vertical unit and composition.

\end{enumerate}

Similarly as before, we can also make $\Gamma$ a simplicial object in the
kind of structures just described, and for the homotopical
part of our discussion we will find it advantageous
to also assume that, in addition with the horizontal category
$Obj_{\Gamma}$  being simplicially constant.

\vspace{3mm}
Now to obtain an analogue of Proposition \ref{pllax*} for lax algebras
enhanced over $\mb$-categories, we will examine the construction leading
up to Proposition \ref{pllax*}, noting along the way how they must
be changed in view of $\mb$-enrichment.

First, we examine $B_{hor}(\Gamma)$. We see that $B_{hor}(Obj_{\Gamma})$
is a $\omc$-algebra over sets, and $B_{hor}(Mor_{\Gamma})$ is a
$\omc$-algebra enhanced over $\mb$ {\em fibered} over $B=B_{hor}(Obj_{\Gamma})$.
We further have associative composition, which is a map of $\omc$-algebras,
$$B_{hor}(Mor_{\Gamma})\boxtimes_{B}B_{hor}(Mor_{\Gamma})\r
B_{hor}(Mor_{\Gamma})$$
(the symbol $\boxtimes_{B}$ indicates $\boxtimes$ applied fiber-wise).
Again, what we want is to be able to extend the product in a way so that
we can replace $\boxtimes_{B}$ with $\boxtimes$. 
Again, we have a category of $\mb$-enhanced $\omc$-algebras fibered over $B$,
and a map of monads $D\r E$ in this category,
$$DX=\cform{\coprod}{n\geq 0}{}\underbrace{X\boxtimes_{B}...
\boxtimes_{B}X}_{\text{$n$ times}},$$
$$EX=\cform{\coprod}{n\geq 0}{}\underbrace{X\boxtimes...
\boxtimes X}_{\text{$n$ times}},$$
may then form the two-sided bar construction of monads \rref{ellax7} and note
that \rref{ellax7} is a $\mb$-enhanced $\omc\Box\mc_{1}$-algebra.
Therefore, we must again find conditions when
the second map \rref{ellax7a} is an equivalence, when $X=B_{hor}(Mor_{\Gamma})$.
The conditions we arrive at are again \rref{ellax8} and \rref{ellax8b}, although
while \rref{ellax8} does not change, in \rref{ellax8b}
the $Hom_{vert}$ sets now denote objects of $s\mb$. When these conditions
are satisfied, we say, again, that $\Gamma$ is {\em distinguished}. We therefore
have an enhanced analogue of Proposition \ref{pllax*}:

\begin{proposition}
\label{pllax**}
Let $\Gamma$ be a distinguished lax $\mc$-algebra enhanced over $\mb$-categories.
Then for $M,N\in Obj_{Obj(\Gamma)}$, $Hom_{vert}(M,N)$ are all naturally equivalent,
and naturally equivalent
to a $\omc\Box\mc_{1}$-algebra enhanced over $\mb$.
\end{proposition}

\vspace{5mm}
\section{Proof of Theorem \ref{t2}}
\label{s4a}

In view of Proposition \ref{pllax**}, it suffices to produce a distinguished lax
$\mc$-algebra $\Gamma$ enhanced in $\mb$-categories, where for $M,N\in Obj_{Obj(\Gamma)}$,
$$Hom_{vert}(M,N)$$
is naturally equivalent to \rref{epr*}.

\vspace{3mm}
Let $Q$ be a cofibrant
fibrant special operad in $s\mb$ fibered over $\mc$, a
cofibrant
operad in simplicial sets with 
$\mc(0)=*$,
$\mc(1)\simeq *$. Define the category $\Gamma$ enhanced in
$s\mb$-categories as follows: If $\mb=K$-modules,
the horizontal category $Obj_{\Gamma}$
has as objects all fibration equivalences of $Q_{1}$-modules
\begin{equation}
\label{eex+}
M\r Q(0)
\end{equation}
where $M$ is cofibrant in the
category of $Q_1$-modules (note that $Q(0)$ is fibrant). If
$\mb=Sets$, drop the requirement that \rref{eex+} be
a fibration. The (horizontal)
morphisms in $Obj_{\Gamma}$ are commutative diagrams
\begin{equation}
\label{eex++}
{\diagram
M\rto\dto_{f}&Q(0)\\
N\urto&
\enddiagram}
\end{equation}
where the map $f:M\r N$ is a cofibration in the category
of $Q_1$-modules.
The vertical $Hom$-sets in $Obj(\Gamma)$ (i.e. between two objects of
the form \rref{eex+})
are
$$Hom_{Q_{1}}(M,M^{\prime}),$$
if $\mb=K$-modules. 

When $\mb=Sets$, we could actually use the {\em space}
\begin{equation}
\label{epri}
Map_{|Q_{1}|}(|M|,|M^{\prime}|).
\end{equation}
However, since we want to stay in the
category of simplicial sets, we get back
by applying the singular set functor to \rref{epri}.

The vertical $Hom$-sets in $Mor(\Gamma)$
are groups of pairs (under addition) of
$$(u,v)\in Hom_{Q_{1}}(M,N)\prod Hom_{Q_{1}}(M^{\prime},N^{\prime})$$
for $K$-modules and
$$(u,v)\in Map_{|Q_{1}|}(|M|,|N|)\prod Map_{|Q_{1}|}(|M^{\prime}|,|N^{\prime}|)
$$
for sets
which make commutative the diagram
$$
\diagram
M\rto^{f}\dto_{u}&M^{\prime}\dto^{v}\\
N\rto_{g}&N^{\prime}.
\enddiagram
$$
Note that the requirement that $f$ be a cofibration implies 
(for $K$-modules) that the
map induced by $f$
$$Hom(M^{\prime},N^{\prime})\r Hom(M,N^{\prime})$$
is a fibration equivalence. Therefore
$$\diagram
Hom(f,g)=Hom(M,N)\protect\cform{\prod}{Hom(M,N^{\prime})}{}
Hom(M^{\prime},N^{\prime})\dto_{\simeq}\\
Hom(M,N)
\enddiagram
$$
is a fibration equivalence. Now we have a diagram (where the arrows
are induced maps)
$$
\diagram
Hom(f,g)\rto\dto_{\simeq}&Hom(M^{\prime},N^{\prime})\dto^{\simeq}\\
Hom(M,N)\rto_{\simeq}& Hom(M,N^{\prime}).
\enddiagram
$$
The bottom row is an equivalence because $N\r N^{\prime}$ is an equivalence,
and by our hypothesis, $N,N^{\prime}$ are fibrant. Similarly for sets.
This already proves condition \rref{ellax8b}.

Now $Q$ being special implies that $\Gamma$ is a lax $\mc$-algebra
enhanced over $s\mb$-categories. In effect, the lax
$\mc$-algebra structure is defined as follows: for $x\in \mc(n)$, and
$Q_{1}$-modules $M_{1},...,M_{n}$, the $x$-product of $M_{1}$,...,$M_{n}$
is
\begin{equation}
\label{eprod+}
Q(n)_{x}\boxtimes_{Q_{1}^{\boxtimes n}}(M_{1}\boxtimes...\boxtimes M_{n}).
\end{equation}
If $Q$ is special, cofibrant and fibrant, the canonical map from \rref{eprod+}
to $Q(0)$ is an equivalence (see \rref{esp}). On morphisms, \rref{eprod+}
preserves $Q_{1}$-cofibrations if $Q$ is cofibrant. For
$\mb=K$-modules, \rref{eprod+} preserves fibrations, because fibrations
of simplicial $K$-modules are precisely onto maps.
Note that this operation is functorial, and on morphisms carries cofibrations
to cofibrations, with $s\mb$ as above.

Thus, the proof of the statement that $\Gamma$ is distinguished (and hence of
Theorem \ref{t2}) is reduced to the following

\begin{proposition}
\label{pspec}
\begin{equation}
\label{epspec1}
BObj_{\Gamma}\simeq *.
\end{equation}
\end{proposition}

\vspace{3mm}
\noindent
{\bf Proof:}
To show \rref{epspec1}, one chooses a particular object
$$M_{0}\r Q(0).$$
For any object of $Obj_{\Gamma}$
$$M\r Q(0),$$
one can obviously choose an arrow $N_{M}\r Q(0)$ in $Obj_{\Gamma}$
together with a diagram in $Obj_{\Gamma}$
$$\diagram
M\rto\drto_{\simeq}&N_{M}\dto^{\simeq}&M_{0}\lto\dlto^{\simeq}\\
&Q(0).&
\enddiagram
$$
Similarly, for a (horizontal) morphism in $Obj_{\Gamma}$
$$\diagram
M\rrto^{f}\drto_{\simeq}&&M^{\prime}\dlto^{\simeq}\\
&Q(0),&
\enddiagram$$
there is an $N_{f}\r Q(0)$ in $Obj_{\Gamma}$ such that, over $Q(0)$, we have
a commutative diagram
$$\diagram
M\rto\drto\ddto&N_{M}\dto&M_{0}\lto\dlto\ddto^{=}\\
\protect{}&N_{f}&\protect{}\\
M^{\prime}\rto\urto&N_{M^{\prime}}\uto&M_{0}.\lto\ulto
\enddiagram
$$
To formalize the procedure this will give, recall the {\em barycentric
subdivision} $C^{\prime}$ of a category $C$ (in our case, $C=Obj_{\Gamma}$):
The category $C^{\prime}$ is a partially ordered
set whose objects are $n$-tuples of
composable arrows in $C$
\begin{equation}
\label{ebar}
{\diagram
x_{0}\rto^{f_{1}}& x_{1}\rto&...\rto^{f_{n}}&x_{n},
\enddiagram}
\end{equation}
and a $k$-tuple of composable morphisms 
$$
\diagram
\cdot\rto^{g_{1}}& \cdot\rto&...\rto^{g_{k}}&\cdot
\enddiagram
$$
is said to be $\leq$ \rref{ebar}
if the $g_{i}$'s are obtained by consecutive compositions 
(or omissions from the beginning or end) of the
$f_{i}$'s. Then there is a canonical functor
$$\Psi:C^{\prime}\r C$$
which on objects is given by
$$(f_{n},...,f_{1})\mapsto Tf_{n}$$
(and by the obvious formula on morphisms). The functor $\Psi$ induces an
equivalence upon applying $B$ (the bar construction).

\vspace{3mm}
Now repeating the procedure
we described constructs a functor 
$$N:(Obj_{\Gamma})^{\prime}\r Obj_{\Gamma},$$
together with natural transformations
\begin{equation}
\label{egamma}
\Psi\r N,
\end{equation}
\begin{equation}
\label{egn}
G\r N
\end{equation}
where $G$ is the constant functor with value in $M_{0}\r Q(0)$, thus showing
that $BObj_{\Gamma}$ is contractible, as claimed.
\qed

\vspace{5mm}

\section{The special property for $k$-cubes}
\label{sspec}

The purpose of this section is to prove Theorem \ref{tspec}.
In this section $\mc_{k}$ will stand for the topological version of
the little $k$-cube operad, and we will also consider the topological
version of the construction \rref{E+}, $\mc=\mc_{k}$. Recall that
$D_{0}(X)=C_{k}X$. It clearly suffices to prove that the map
\begin{equation}
\label{E++}
\phi:B(D_{n}X,(D_{1}X)^{\times n},(D_{0}X)^{\times n})
\r C_{k}X\times B(\mc_{k}(n), \mc_{k}(1)^{\times n},\mc_{k}(0)^{\times n})
\end{equation}
is an equivalence, where the first coordinate of the map $\phi$ is
given by composition, and the second map by the forgetful map
$$D_{\ell}(X)\r \mc_{k}(\ell).$$
In effect, to get from spaces to simplicial sets, we may apply the
singular set functor, and to get to $K$-modules, we may further 
apply the free $K$-module functor.

Now since the left hand side of \rref{E++} obviously preserves weak 
equivalences, we can further replace the terms of \rref{E++} as follows.
First, let $\mc^{\prime}_{k}(\ell)\subset \mc_{k}(\ell)$ consist of all
$\ell$-tuples of little cubes
$\alpha_{1},...\alpha_{\ell}:I^{k}\r I^{k}$
such that
$Im(\alpha_{i})\subset Int(I^{k})$ and the images $Im(\alpha_{i})$ are
disjoint. Now let $M(\ell)$ be the space of pairs $(\alpha,m)$ where
$\alpha=(\alpha_{1},...,\alpha_{\ell})\in\mc_{k}(\ell)$ and $m$ is a set
of unordered $X$-decorated points in $I^{k}-
\protect\cform{\bigcup}{i=1}{\ell}
Im(\alpha_{i})$ (with the usual configuration space topology).
Let $M^{\prime}(\ell)$ be defined in the same way as $M(\ell)$, with
$\mc_{k}(\ell)$ replaced by $\mcpk(\ell)$:
$$M^{\prime}(\ell)=M(\ell)\times_{\mc_{k}(\ell)}\mcpk(\ell).$$
Then we have equivalences 
$$D_{n}(X)\r M(n) \leftarrow M^{\prime}(n)$$
(the first map replaces a little cube decorated by an element of $X$
by its center). Thus, we can further restate the claim that \rref{E++}
is an equivalence as follows:

\begin{proposition}
\label{pspec1}
The map
\begin{equation}
\label{Ei}
\kappa:B(M^{\prime}(\ell), M(1)^{\times \ell}, M(0)^{\times \ell})
\r M(0)\times B(\mcpk(\ell),\mc_{k}(1)^{\times \ell},\mc_{k}(0)^{\times
\ell})
\end{equation}
where the first map is by composition and the second map is by projection
$M\r\mc$, is an equivalence.
\end{proposition}

\vspace{3mm}
\noindent
{\bf Remark:}
Recall that $\mcpk(\ell)\simeq \mc_{k}(\ell)$, $\mc_{k}(1)\simeq *$,
so the right hand side of \rref{Ei} is weakly equivalent to $M(0)\times
\mc_{k}(\ell)$.

\vspace{3mm}
\noindent
{\bf Proof of Proposition \ref{pspec1}:}
First define $N(\ell)$ as the pullback
\begin{equation}
\label{EE+}
N(\ell)=(\mcpk(\ell)\times\mc_{k}(1)^{\ell})\times_{\mcpk(\ell)}M^{\prime}
(\ell)
\end{equation}
where the map
$$\mcpk(\ell)\times \mc_{k}(1)^{\ell}\r \mcpk(\ell)$$
is by composition. Then $N(\ell)$ enjoys a right $M(1)^{\ell}$-action
where $M(1)^{\ell}$ acts by composition on $M^{\prime}(\ell)$ and
by the forgetful map together with internal composition on $\mcpk(\ell)
\times \mc_{k}(1)^{\ell}$:
$$
\diagram
(\mcpk(\ell)\times\mc_{k}(1)^{\ell})\times M(1)^{\ell}\rto
&\mcpk(\ell)\times \mc_{k}(1)^{\ell}\times \mc_{k}(1)^{\ell}\dto^{Id\times
\gamma^{\ell}}\\ \protect{}&\mcpk(\ell)\times\mc_{k}(1)^{\ell}.
\enddiagram
$$
Furthermore, the projection 
$$N(\ell)\r M^{\prime}(\ell)$$
is obviously an equivalence (by contracting $\mc_{k}(1)^{\ell}$),
so we may replace $M^{\prime}(\ell)$ by $N(\ell)$ in the statement of
the Proposition.

\vspace{3mm}
Now filter $N(\ell)$ by closed subspaces $N_{q}(\ell)$ consisting of
all triples
\begin{equation}
\label{EE++}
(\alpha,\beta,x),
\end{equation}
$\alpha=(\alpha_{1},...,\alpha_{\ell})\in\mcpk(\ell)$,
$\beta=(\beta_{1},...,\beta_{\ell})\in \mcpk(1)^{\ell}$,
$x\in M^{\prime}(\ell)$
where the number of $X$-decorated points of $x$ contained in
$$I^{k}-\cform{\bigcup}{i=1}{\ell}Im(\alpha_{i})$$
is $\leq q$.

\begin{lemma}
\label{lspec1}
Suppose $A\subset X$ is an $M$-equivariant weak NDR pair with
Urysohn function and homotopy $(u,h)$ where $M$ is a monoid acting
on the right, such that the following conditions are satisfied:
\begin{enumerate}
\item
\label{ii1}
$X-A=V\times M$ for some $V\subset X-A$ closed.

\item
\label{ii2}
$h_{t}(A\cup V)\subset A\cup V$ for all $t$. (Note that $A\cup V
=A\cup Cl(V)$. $Cl$ denotes closure.)

\end{enumerate}
Then there is a natural $M$-equivariant homotopy equivalence
$$Cofiber(A\subset X)\r ((A\cup V)/A)\wedge M_{+}.$$
\end{lemma}

\Proof
Because of the weak NDR property, we have
$$Cofiber(A\subset X)\simeq X/A,$$
so it suffices to prove
\begin{equation}
\label{EL+}
X/A\simeq (A\cup V)/A\wedge M_{+}.
\end{equation}
To get a map $\leftarrow$, note that the inclusion
$$A\cup V\r X$$
extends to an equivariant map
$$(A\cup V)\times M\r X\r X/A$$
which clearly annihilates $A\times M$, thus inducing
$$\psi: (A\cup V)\times M/A\times M\r X/A.$$
To get a map $\r$ in \rref{EL+},
use the map $\phi$ induced by $h_{1}$: to see that it is continuous,
let
$$U=\{x\in X| u(x)<1\}.$$
Then $\phi$ is constant (hence continuous) on $U$, but also continuous
on $X-A$. Now $\{U,X-A\}$ is an open cover of $X$.

Now since $\psi$ is obviously a bijection, we can define both homotopies
$\phi\psi\simeq Id$, $\psi\phi\simeq Id$ as $h_{t}$. Then $\psi\phi$ is
a quotient of $h_{t}$ (with topology), and thus is continuous.
To see that 
\begin{equation}
\label{Eu1}
h_{t}:\phi\psi\simeq Id 
\end{equation}
is continuous, note that we have
a continuous map
$$h_{t}:A\cup V\r A\cup V,$$
and hence
\begin{equation}
\label{Eu}
unit\circ h_{t}:A\cup V\r (A\cup V)\wedge M_{+}.
\end{equation}
Then \rref{Eu1} is the free extension of \rref{Eu}.
\qed

Now below we shall prove

\begin{lemma}
\label{lspec2}
The pair $N_{q-1}(\ell)\subset N_{q}(\ell)$, with $M=M(1)^{\ell}$, satisfy
the hypotheses of Lemma \ref{lspec1}.
\end{lemma}

Thus, so do $B(N_{q-1}(\ell),M(1)^{\ell},M(1)^{\ell})
\subset B(N_{q}(\ell),M(1)^{\ell},M(1)^{\ell}).$

Now consider the diagram
\begin{equation}
\label{EE1}
{\diagram
B(N_{q-1}(\ell),M(1)^{\ell},M(0)^{\ell})\rto\dto
& B(N_{q}(\ell),M(1)^{\ell},M(0)^{\ell})\dto\\
N_{q-1}(\ell)\times_{M(1)^{\ell}}M(0){\ell}\rto
& N_{q}(\ell)\times_{M(1)^{\ell}}M(0)^{\ell}.
\enddiagram}
\end{equation}
By Lemma \ref{lspec1}, the cofibers of both rows are naturally equivalent
to
$$((A\cup V)/A)\wedge M(0)^{\ell}_{+}.$$
Thus, if the left column of \rref{EE1} is an equivalence, so is the right
column. Thus, inductively,
$$B(N(\ell),M(1)^{\ell},M(0)^{\ell})
\simeq N(\ell)\times_{M(1)^{\ell}}M(0)^{\ell}\cong \mcpk(\ell)\times M(0).$$
The equivalences are easily checked to be compatible with the required maps.
\qed

\vspace{3mm}
\noindent
{\bf Proof of Lemma \ref{lspec2}:}
To simplify notation, we shall identify little cubes with their images.
For a little cube $\alpha$ in $I^{k}$, and for $t\in \R_{>0}$, let
$t\alpha$ be the cube with the same center which is $t$ times larger.

Now first note that there is a continuous function
$$\lambda:\mcpk(\ell)\r \R_{>1}$$
such that for $\alpha=(\alpha_{1},...,\alpha_{\ell})\in\mcpk(\ell)$,
$\lambda=\lambda(\alpha)$,
\begin{equation}
\label{EE2}
(\lambda\alpha_{1},...,\lambda\alpha_{\ell})\in\mcpk(\ell).
\end{equation}
To this end, for every $\beta\in\mcpk(\ell)$, there is an open
neighbourhood $U$ of $\beta$ and a {\em constant}
$\mu>1$
which works as $\lambda$ in \rref{EE2} for $\alpha\in U$.
Thus, since $\mcpk(\ell)$ is paracompact, the function $\lambda$ can 
constructed by partition of unity.

\vspace{3mm}
Now also note that for two pairs of cubes $\alpha\supset\beta$,
$\gamma\supset\delta$ (not equal), there is a canonical homeomorphism
$$\Phi^{\alpha\beta}_{\gamma\delta}:Cl(\beta-\alpha)\r
Cl(\delta-\gamma).$$
First, there are canonical linear homeomorphisms identifying the boundaries
$\partial\alpha,\partial\beta,\partial\gamma,\partial\delta$. Connecting
two corresponding points on $\partial\alpha$, $\partial\beta$ creates a line
segment; map this segment linearly onto the line 
segment obtained by identifying
the corresponding points on $\partial\gamma$, $\partial\delta$.
Let also $(\Phi^{\alpha\beta}_{\gamma\delta})_{*}$ be the map induced by
$\Phi^{\alpha\beta}_{\gamma\delta}$ on configuration spaces.

Then define the homotopy $h$ in $N_{q}(\ell)$ (see \rref{EE+}) by
$$h_{t}(\alpha,1,m)=
((\alpha_{i})_{i=1}^{\ell},((1-t/3))_{i=1}^{\ell},
((\Phi^{\lambda(\alpha)\alpha_{i},\alpha_{i}}_{
\lambda(\alpha)\alpha_{i},(1-t/3)\alpha_{i}})_{i=1}^{\ell})_{*}m),$$
and extend $M(1)^{\ell}$-equivariantly. We can let 
$$V$$
consist of all triples $(\alpha,1,m)$ (see \rref{EE+})
where $m$ contains exactly $q$
$X$-decorated points. Define
$$
u(\alpha,\beta,m)=sup\left\{t\left|
\parbox{3in}{there are exactly $q$ $X$-decorated
points in 
$((\Phi^{\lambda(\alpha)\alpha_{i},\alpha_{i}}_{
\lambda(\alpha)\alpha_{i},(1-t/3)\alpha_{i}})_{i=1}^{\ell})_{*}m$
which are in $I^{k}-\cform{\bigcup}{i=1}{\ell}\alpha_{i}$.
}\right.\right\}
$$
\qed

\vspace{5mm}
{\sc
Po Hu

Department of Mathematics

University of Chicago

5734 S University Ave.

Chicago, IL 60637

\vspace{5mm}
Igor Kriz

Department of Mathematics 

University of Michigan

2074 E Hall, 525 E University Ave.

Ann Arbor, MI 48109-1109

\vspace{5mm}
Alexander  A. Voronov

Department of Mathematics

University of Minnesota

127 Vincent, 206 Church St. SE

Minneapolis, MN 55455-0487

}
\end{document}